\documentclass[12pt]{amsart}

\usepackage{enumerate}
\usepackage{latexsym}
\usepackage[pdftex]{graphicx}
\usepackage{amssymb}
\usepackage{amsmath}
\usepackage{xcolor}
\usepackage{stackrel}
\usepackage[normalem]{ulem}
\usepackage{thumbpdf}
\usepackage[pagebackref]{hyperref} 

\newtheorem*{theorem}{Theorem}
\newtheorem{lemma}{Lemma}
\newtheorem{proposition}{Proposition}
\newtheorem{remark}{Remark}
\newtheorem{example}{Example}

\newtheorem{definition}{Definition}
\newtheorem{corollary}{Corollary}

\newtheorem{problem}{Problem}
\newtheorem*{problem*}{Problem}

\newtheorem{innercustomthm}{Problem}
\newenvironment{customthm}[1]
  {\renewcommand\theinnercustomthm{#1}\innercustomthm}
  {\endinnercustomthm}

\newcommand{\beq}{\begin{equation}}
\newcommand{\eeq}{\end{equation}}
\newcommand{\beqna}{\begin{eqnarray*}}
\newcommand{\eeqna}{\end{eqnarray*}}
\newcommand{\beqn}{\begin{equation*}}
\newcommand{\eeqn}{\end{equation*}}
\newcommand{\bp}{\begin{proof}}
\newcommand{\ep}{\end{proof}}
\newcommand{\bprop}{\begin{proposition}}
\newcommand{\eprop}{\end{proposition}}
\newcommand{\bt}{\begin{theorem}}
\newcommand{\et}{\end{theorem}}
\newcommand{\bex}{\begin{example}}
\newcommand{\eex}{\end{example}}
\newcommand{\bc}{\begin{corollary}}
\newcommand{\ec}{\end{corollary}}
\newcommand{\bl}{\begin{lemma}}
\newcommand{\el}{\end{lemma}}
\newcommand{\bprob}{\begin{problem}}
\newcommand{\eprob}{\end{problem}}
\newcommand{\br}{\begin{remark}}
\newcommand{\er}{\end{remark}}
\newcommand{\bd}{\begin{definition}}
\newcommand{\ed}{\end{definition}}

\begin{document}

\title
[Fifth and eighth Busemann-Petty problems near the ball]
{A  solution to the fifth and the eighth Busemann-Petty problems in a small neighborhood of the Euclidean ball}

\author[M.A. Alfonseca]{M. Angeles Alfonseca}
\address{Department of Mathematics, North Dakota State University\\
Fargo, ND 58108, USA} \email{maria.alfonseca@ndsu.edu}

\author[F. Nazarov]{Fedor Nazarov}
\address{Department of Mathematical Sciences, Kent State University,
Kent, OH 44242, USA} \email{nazarov@math.kent.edu}

\author[D. Ryabogin]{Dmitry Ryabogin}
\address{Department of Mathematical Sciences, Kent State University,
Kent, OH 44242, USA} \email{ryabogin@math.kent.edu}

\author[V. Yaskin]{Vladyslav Yaskin}
\address{Department of Mathematical and Statistical Sciences, University of Alberta,
Edmonton, Canada} \email{vladyaskin@math.ualberta.ca}

\thanks{The first author is supported in part by the Simons Foundation Grant 711907.
	The second  and the third authors are  supported in
part by U.S.~National Science Foundation Grants  DMS-1900008 and DMS-1600753. The fourth author is supported by NSERC}


\begin{abstract}
We show that the fifth and the eighth Busemann-Petty  problems  have positive solutions for bodies that are sufficiently close to the Euclidean ball in the Banach-Mazur distance.
\end{abstract}

\maketitle

\section{Introduction}

In  1956, Busemann and Petty  \cite{BP} posed ten problems
about symmetric convex bodies, of which only the first one has been
solved (see \cite{K}). Their fifth and eighth problems are as follows.

\begin{customthm}{5}\label{bp5}
 
 If for an origin-symmetric convex body $K\subset{\mathbb R^n}$, $n \geq 3$, we have
\begin{equation}\label{as1}
h_K(\theta) \textnormal{vol}_{n-1}(K\cap\theta^{\perp})=C\qquad\forall\theta\in S^{n-1},
\end{equation}
where the constant $C$ is independent of $\theta$,  must $K$ be an ellipsoid?

\end{customthm}

Here $S^{n-1}$ is the unit sphere in $\mathbb R^n$,
 $\theta^{\perp}=\{x\in{\mathbb R^n}:\langle
  x,\theta\rangle =0\}$ is the hyperplane orthogonal to the 
direction $\theta\in S^{n-1}$ and passing through the origin, and
  $h_K(\theta)=\max\limits_{x\in K} \langle x ,  \theta\rangle$ is the support function of the convex body $K\subset {\mathbb R^n}$.

\begin{customthm}{8}\label{bp8}
 
 If for an origin-symmetric convex body $K\subset{\mathbb R^n}$, $n \geq 3$, we have
\begin{equation}\label{as11}
f_K(\theta)=C(\textnormal{vol}_{n-1}(K\cap \theta^{\perp}))^{n+1}\qquad\forall \theta\in S^{n-1},
\end{equation}
where the constant $C$ is independent of $\theta$,  must $K$ be an ellipsoid?

\end{customthm}

Here $f_K$ is the curvature function of $K$, which is the reciprocal of the Gaussian curvature viewed as a function of the unit normal vector (see \cite[pg. 419]{Sch}).

The Euclidean ball clearly satisfies \eqref{as1} and \eqref{as11}. If a body $K$ satisfies \eqref{as1},  then  so does $TK$ for any linear transformation $T\in GL(n)$ (with  constant $C \cdot |\det T|$). Similarly,  if a body $K$ satisfies  \eqref{as11}, then  so does $TK$ for any linear transformation $T\in GL(n)$ (with  constant $C \cdot |\det T|^{1-n}$). Hence, \eqref{as1} and \eqref{as11} are satisfied by ellipsoids.

\bigskip

In this paper we prove the following result.  

\bt\label{mtH1}
Let $n\ge 3$. If  an origin-symmetric convex body $K\subset {\mathbb R^n}$ satisfies \eqref{as1} or \eqref{as11} and is sufficiently close to the  Euclidean ball  in the Banach-Mazur metric, then $K$ must be an ellipsoid.
\et

In dimension 2, there are  convex bodies  satisfying \eqref{as1} that are not ellipses but, nevertheless,  can be arbitrarily close to the unit disc. The curve bounding such a body is  a so-called {\it Radon curve}, see \cite{D}. On the other hand,  the only convex bodies satisfying \eqref{as11} in dimension 2 are the ellipses  \cite[Theorem 5.6]{P}.

\section{Invariance of Busemann-Petty problems under linear transformations}

Both Busemann-Petty problems are invariant under linear transformations in the sense that if a symmetric convex body $K$ satisfies \eqref{as1} or \eqref{as11}, then so does $TK$ where $T$ is an invertible linear map from ${\mathbb R^n}$ to itself.

This statement is almost obvious for Problem \ref{bp5}. Indeed, let $H$ be any hyperplane in ${\mathbb R^n}$ passing through the origin and let $H_s$ be a support hyperplane of $K$ parallel to $H$. Consider any point $p\in K\cap H_s$ and the cone $C_{K,H}$ with the base $K\cap H$ and the vertex $p$. Note that due to the symmetry of $K$ and the fact that $H_s$ is parallel to $H$, the volume $\textrm{vol}_n(C_{K,H})$ of this cone is independent of the particular choice of $H_s$ and $p$. Moreover, we clearly have $C_{TK, TH}=T(C_{K,H})$, so $\textrm{vol}_n(C_{TK, TH})=|\det T| \textrm{vol}_n(C_{K,H})$.  Since for $H=\theta^{\perp}$ this volume can be expressed as $\textrm{vol}_n(C_{K,H})=\frac{1}{n}\textrm{vol}_{n-1}(K\cap \theta^{\perp})h_K(\theta)$,  we see that property \eqref{as1} is merely the statement that 
$\textrm{vol}_n(C_{K,H})$ is independent of the choice of the hyperplane $H$ (this was exactly how the fifth Busemann-Petty problem was originally formulated in \cite{BP}).

The invariance of \eqref{as11} under linear transformations is somewhat less transparent.  When $K$ has smooth non-degenerate  $C^2$-boundary with strictly positive Gaussian curvature at each point, we can restate it as follows. 

Let, as before, $H$ be an arbitrary hyperplane passing through the origin, let $H_s$ be one of the two supporting hyperplanes of $K$ parallel to $H$, and let $p\in K\cap H_s$. For $t\in (0,1)$, let $H^t$ be the hyperplane between $H$ and $H_s$ parallel to $H$ such that the distance between $H_s$ and $H^t$ is $t$ times the distance $d$ between $H_s$ and $H$. Then, for small $t$, the $(n-1)$-dimensional volume $\textrm{vol}_{n-1}(K\cap H^t)$ is approximately proportial to $\frac{  t^{\frac{n-1}{2}  } d^{\frac{n-1}{2}} }{\sqrt{G(p)}}$ where $G(p)$ is the Gaussian curvature of $\partial K$ at $p$.

Note now that $\frac{\textrm{vol}_{n-1}(K\cap H^t)}{\textrm{vol}_{n-1}(K\cap H)}$ is invariant under linear transformations and
$d\,\textrm{vol}_{n-1}(K\cap H)=n\,\textrm{vol}_n(C_{K,H})$ is multiplied by $|\det T|$ when we replace $K$ by $TK$ and $H$ by $TH$. Thus, $G(p)\textrm{vol}_{n-1}(K\cap H)^{n+1}$
equals (up to a constant factor depending on the dimension $n$ only) 
\[
\lim\limits_{t\to 0}t^{n-1}\textrm{vol}_n(C_{K,H})^{n-1}\Big[\frac{\textrm{vol}_{n-1}(K\cap H)}{\textrm{vol}_{n-1}(K\cap H^t)} \Big]^2,
\]
and thus is multiplied by $|\det T|^{n-1}$ when we replace $K$ by $TK$ and $H$ by $TH$.

In general, it is unclear to us what degree of smoothness Busemann and Petty assumed when posing Problem \ref{bp8}. We will handle the most general case, when \eqref{as11} is understood in the sense that the surface area measure of $K$ is absolutely continuous with respect to the $(n-1)$-dimensional Lebesgue measure on the unit sphere and its Radon-Nikodym density is equal to the right-hand side. In this case, 
the geometric meaning of \eqref{as11} is less transparent but the invariance of \eqref{as11} under linear transformations still follows from the computations in  \cite{Lu}. The reader will lose almost nothing, however, by assuming that $\partial K$ is smooth and non-degenerate, but $K$ is close to the unit ball only in the Banach-Mazur distance and not in $C^2$.

\section{From the Banach-Mazur distance to the Hausdorff one}\label{BMF}

Applying an appropriate linear transformation, we can assume that the constants in \eqref{as1} and \eqref{as11} are equal to those for the unit ball $B^n_2$ and that $(1-\varepsilon)r B^n_2\subset K\subset (1+\varepsilon)rB^n_2$ for some $r>0$ with some small $\varepsilon>0$.

Our task here will be to show that $r$ must be close to $1$, i.e., $K$ must be close to the unit Euclidean ball $B^n_2$ in the Hausdorff metric. We have
\begin{equation}\label{votvot1}
(1-\varepsilon)rh_{B^n_2}\le h_K\le (1+\varepsilon)rh_{B^n_2}
\end{equation}
and
\begin{equation}\label{votvot2}
(1-\varepsilon)^{n-1}r^{n-1}\textrm{vol}_{n-1}(B^n_2\cap\theta^{\perp})\le    \textrm{vol}_{n-1}(K\cap\theta^{\perp})\le 
\end{equation}
\[
(1+\varepsilon)^{n-1}r^{n-1}\textrm{vol}_{n-1}(B^n_2\cap\theta^{\perp}).
\]
In the case of \eqref{as1}, combining \eqref{votvot1} and \eqref{votvot2}  with the equation 
\[
h_K(\theta)\textrm{vol}_{n-1}(K\cap\theta^{\perp})=h_{B^n_2}(\theta)\textrm{vol}_{n-1}(B^n_2\cap\theta^{\perp}),
\]
we obtain $(1-\varepsilon)^nr^n\le 1\le (1+\varepsilon)^nr^n$, i.e., $\frac{1}{1+\varepsilon}\le r\le \frac{1}{1-\varepsilon}$.

\bigskip

In the case of \eqref{as11}, we can integrate both sides with respect to the $(n-1)$-dimensional Lebesgue measure on $S^{n-1}$ to conclude  (see \cite{Sch},  Section 5.3.1)   that
\begin{equation}\label{surface}
\Sigma(K)=\int_{S^{n-1}}f_K(\theta)dm_{n-1}(\theta)=
\end{equation}
$$
c_n\int_{S^{n-1}} \left( \textrm{vol}_{n-1}(K\cap\theta^{\perp})\right)^{n+1}dm_{n-1}(\theta),
$$
where $\Sigma(K)$ is the surface area of $\partial K$ and $c_n$ is defined by 
\[
\Sigma(B^n_2)=c_n\int_{S^{n-1}} \left(\textrm{vol}_{n-1}(B^n_2\cap\theta^{\perp})\right)^{n+1}dm_{n-1}(\theta).
\]
From our assumption $(1-\varepsilon)r B^n_2\subset K\subset (1+\varepsilon) rB^n_2$, it follows that 
\[
(1-\varepsilon)^{n-1}r^{n-1}\Sigma(B^n_2)\le   \Sigma(K)\le (1+\varepsilon)^{n-1}r^{n-1}\Sigma(B^n_2),
\]
which, together with  \eqref{votvot2} and \eqref{surface}, gives
\[
(1-\varepsilon)^{n-1}r^{n-1}\le (1+\varepsilon)^{(n-1)(n+1)}r^{(n-1)(n+1)}
\]
and 
\[
(1+\varepsilon)^{n-1}r^{n-1}\ge (1-\varepsilon)^{(n-1)(n+1)}r^{(n-1)(n+1)},
\]
i.e.,
\[
\frac{1-\varepsilon}{(1+\varepsilon)^{n+1}}\le r^n\le \frac{1+\varepsilon}{(1-\varepsilon)^{n+1}}.
\]

\section{The isotropic position}

We have seen in the previous section that, without loss of generality, we may assume that $(1-\varepsilon)B^n_2\subset K\subset (1+\varepsilon)B^n_2$. However, this requirement still leaves some freedom as to what affine image of $K$ to choose.  In this section we will reduce this freedom even further by putting $K$ into the so-called isotropic position, i.e., the position where
\[
\int_K\langle x,y\rangle^2dy=c\,|x|^2\qquad\qquad \forall x\in {\mathbb R^n}.
\]

The existence of such a position is well known and easy to derive (see \cite{BGVV},  Section 2.3.2). Indeed, for an arbitrary symmetric convex body $K$, the mapping
\[
{\mathbb R^n}\ni x \, \mapsto \int_K\langle x,y\rangle^2 dy=\sum\limits_{i,j} \Big(\int_Ky_iy_jdy\Big)x_ix_j
\]
is a positive-definite quadratic form. Thus, it can be written as $\langle Sx,x\rangle$, where $S$ is a self-adjoint positive definite operator on ${\mathbb R^n}$. 

Moreover, if $K=B^n_2$, then $S=c_n I$ for some $c_n>0$. If  $(1-\varepsilon)B^n_2\subset K$, then
\[
\langle Sx,x\rangle=\int_K\langle x,y\rangle^2dy\ge \int_{(1-\varepsilon)B^n_2}\langle x,y\rangle^2dy=(1-\varepsilon)^{n+2}c_n|x|^2
\]
and, similarly, if $K\subset (1+\varepsilon)B^n_2$, then 
\[
\langle Sx,x\rangle\le (1+\varepsilon)^{n+2}c_n|x|^2.
\]
Thus, setting $\widetilde{S}=c_n^{-1}S$, we have 
\[
(1-\varepsilon)^{n+2}|x|^2  \le \langle \widetilde{S}x,x\rangle\le (1+\varepsilon)^{n+2}|x|^2.
\]
It follows that 
\[
(1-\varepsilon)^{n(n+2)}\leq \det (\widetilde{S})\le (1+\varepsilon)^{n(n+2)}
\]
and 
\[
\|\widetilde{S}\|\le (1+\varepsilon)^{n+2},\quad \|\widetilde{S}^{-1}\|\le (1-\varepsilon)^{-(n+2)}, 
\]
whence the operator $T=\sqrt{\det (\widetilde{S})^{-\frac{1}{n}}\widetilde{S}}$ satisfies 
\[
\det T=1,\qquad\qquad \|T\|,\,\|T^{-1}\|\le \Big(\frac{1+\varepsilon}{1-\varepsilon}\Big)^{\frac{n+2}{2}},
\]
and $T^{-1}ST^{-1}$ is a multiple of the identity. 

The body $\widetilde{K}=T^{-1}K$ satisfies
\[
\int_{\widetilde{K}} \langle x,y\rangle^2dy=\int_{K} \langle x,T^{-1}y\rangle^2dy=\int_{K} \langle T^{-1}x,y\rangle^2dy= 
\]
\[
\langle ST^{-1}x, T^{-1}x\rangle=
\langle T^{-1}ST^{-1}x, x\rangle=c|x|^2
\]
for some $c>0$, while we also have
\[
(1-\varepsilon)\Big(\frac{1-\varepsilon}{1+\varepsilon}\Big)^{\frac{n+2}{2}}B^n_2\subset T^{-1}(1-\varepsilon) B^n_2\subset T^{-1}K
\]
\[
\subset 
T^{-1}(1+\varepsilon) B^n_2\subset (1+\varepsilon)\Big(\frac{1+\varepsilon}{1-\varepsilon}\Big)^{\frac{n+2}{2}}B^n_2.
\]

\section{The analytic reformulation}\label{Fedjaletaet}

Let $\rho_K$, $h_K: S^{n-1}\to{\mathbb R}$ be the radial and the support functions of the convex body $K$ respectively, i.e., 
\[
\rho_K(\theta)=\max\{t>0:\, t\theta\in K   \}
\]
and
\[
h_K(\theta)=\max\{\langle x, \theta\rangle:\, x\in K  \}.
\]

The $(n-1)$-dimensional volume of the section $K\cap\theta^{\perp}$ is given by
\[
\textrm{vol}_{n-1}(K\cap\theta^{\perp})=c_n \mathcal{R}\Big[\rho_K^{n-1}  \Big],
\]
where  $c_n$ is a positive constant depending on the dimension $n$ only and $\mathcal{R}$ is the Radon transform on $S^{n-1}$, i.e., 
\[
\mathcal{R}f(\theta)=\int_{S^{n-1}\cap\theta^{\perp}}f(\xi)d\sigma(\xi)
\]
with $\sigma$ being the $(n-2)$-dimensional Lebesgue measure on $S^{n-1}\cap\theta^{\perp}$ normalized by the condition $\sigma(S^{n-1}\cap\theta^{\perp})=1$, i.e., $\mathcal{R}1=1$. Thus, condition \eqref{as1} can be rewritten as $h_K\mathcal{R}\Big[ \rho_K^{n-1}\Big]=C$,
where, due to the normalization made at the beginning of Section \ref{BMF}, the constant $C$ should be the same as for the unit ball $B^n_2$, i.e., $C=1$. So, we arrive at the equation
\begin{equation}\label{nunu}
h_K=\Big(\mathcal{R}\Big[ \rho_K^{n-1}\Big]\Big)^{-1}.
\end{equation}

Rewriting \eqref{as11} in terms of $h_K$ and $\rho_K$ is trickier. The right-hand side presents no problem: it is just proportional to  $\Big(\mathcal{R}\Big[\rho_K^{n-1}\Big]\Big)^{n+1}$.  So, the equation becomes 
$f_K=C\Big(\mathcal{R}\Big[\rho_K^{n-1}\Big]\Big)^{n+1}$. Due to the normalization made at the beginning of Section \ref{BMF}, the constant $C$ should be the same as for the unit ball $B^n_2$, i.e., $C=1$.
However, $f_K$ can be readily expressed in terms of $h_K$ only if $h_K$ is $C^2$ and we have made no such assumption.

The expression for $f_K$  in the $C^2$-case  can be written as
$f_K=Ah_K$
where the operator $A$ is defined as follows. For a function $h\in C^2(S^{n-1})$ denote by $H(x)$ its degree $1$
homogeneous extension  to the entire space
(i.e., $H(x)=|x|h(\frac{x}{|x|})$ for $x\neq 0$). Let $\widehat{H}=(H_{x_ix_j})_{i,j=1}^n$ be the Hessian of $H$ and let $\widehat{H}_j$ be the matrix obtained from $H$ by removing the $j$-th row and the $j$-th column.
Let $Ah$ be the restriction
of $\sum\limits_{j=1}^n\det \widehat{H}_j$ to the unit sphere $S^{n-1}$
 (see   \cite{Sch}, Corollary 2.5.3).

We shall show that when $\rho_K$ is close to $1$, we can solve the equation 
\[
Ah=\Big(\mathcal{R}\Big[ \rho_K^{n-1}\Big]\Big)^{n+1}
\]
with $h$ close to $1$ in $C^2$. This $h$ will determine a convex body $L$ that satisfies $f_L=f_K$. By the uniqueness theorem  (see \cite{Sch}, Theorem 8.1.1)  we will then conclude that $K=L$, so $h_K=h$ and the smoothness of $h_K$ will be justified a posteriori. Thus, it  will be possible to rewrite \eqref{as11} as
\begin{equation}\label{2'}
Ah_K=\Big(\mathcal{R}\Big[ \rho_K^{n-1}\Big]\Big)^{n+1}.
\end{equation}

\section{Maximal function}

For $e \in S^{n-1}$, $\vartheta\in (0,\pi]$, let $S_\vartheta(e)=\{e'\in S^{n-1},\,\langle e, e'\rangle\ge \cos\vartheta\}$ denote the spherical cap centered at $e$ with spherical radius $\vartheta$. The spherical Hardy-Littlewood maximal function is defined by 
\[
Mf(e)=\max\limits_{\vartheta\in (0,\pi]}\frac{1}{\sigma(S_\vartheta(e))}\int_{S_\vartheta(e)} |f(x)|d\sigma(x),\qquad f\in L^1(S^{n-1}),
\]
where $\sigma$ is the surface measure on $S^{n-1}$ normalized by the condition $\sigma(S^{n-1})$ $=1$. It is well known that
$M$ is bounded as an operator from $L^2(S^{n-1})$ to itself  (see \cite{Kn}).

\begin{figure}[h]
	\centering
	\includegraphics[height=4.0in]{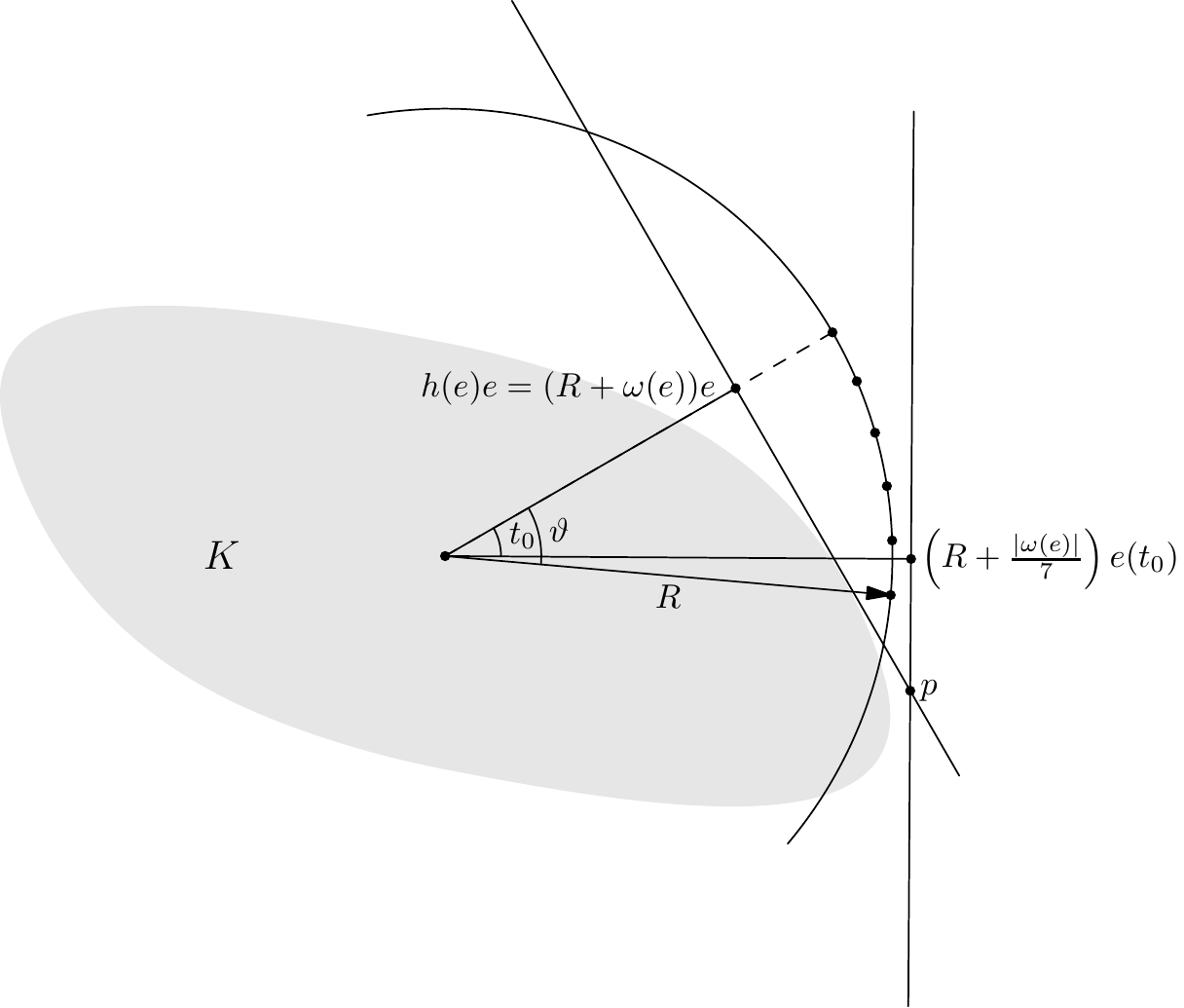} 
	\caption{The body $K$, the lines $\langle x,e\rangle=R+\omega(e)$, $\langle x, e'(t_0)\rangle =R+\frac{1}{7}|\omega(e)|$, and the point $p$}
	\label{bodyK}
\end{figure}

\bl\label{Fl5}
Let $K$ be a $2$-dimensional origin-symmetric convex body  and let $R$ be a positive real number. Let $h_K=R+\omega$ be the support function of $K$ and let $e\in S^1$ be a unit vector. Assume that $h_K(e)\le R\cos\vartheta$ for some $\vartheta\in (0,\frac{\pi}{2})$. Denote by $e'(t)$ the unit vector situated clockwise from $e$ and making an angle $t$ with $e$. Then
\[
|\omega(e)|\le \frac{35}{\vartheta}\int_{\frac{\vartheta}{5}}^\vartheta|\omega(e'(t))|dt.
\]
\el
\bp
Note that the hypothesis $h_K(e)\le R\cos\vartheta$ implies that $\omega(e)<0$.
If $\omega(e'(t))\ge \frac{1}{7}|\omega(e)|$ for all $t\in [\frac{4\vartheta}{5},\vartheta]$, the inequality  clearly holds. Assume now that $\omega(e'(t_0))<\frac{1}{7}|\omega(e)|$ for some $t_0\in [\frac{4\vartheta}{5},\vartheta]$. Let $p$ be the intersection point of the lines $\langle x,e\rangle=R+\omega(e)$ and $\langle x, e'(t_0)\rangle =R+\frac{1}{7}|\omega(e)|$. Then $p$ lies clockwise from $e$ and, since $|p|>R$ and $\langle p,e\rangle=h_K(e)<R\cos\vartheta$, the angle $\alpha$ between $p$ and $e$ is at least $\vartheta$ (see Figure \ref{bodyK}). Also, since $\langle p,e\rangle=h_K(e)>0$, we have $\alpha<\frac{\pi}{2}$.

Since $K$ is contained entirely in the angle $\langle x, e\rangle\le R+\omega(e)$, $\langle x, e'(t_0)\rangle\le R+\frac{1}{7}|\omega(e)|$ with vertex $p$, we have $h_K(e'(t))\le \langle p,e'(t)\rangle=|p|\cos(\alpha-t)$ for all $t\in [0,t_0]$.

We shall now use the following elementary property of the cosine function: if $\gamma$, $\delta$$>0$ and $[\gamma-\delta,\gamma+\delta]\subset [0,\frac{\pi}{2}]$, then $\cos\beta\le \frac{3\cos(\gamma-\delta)+\cos(\gamma+\delta)}{4}$ for all $\beta\in [\gamma,\gamma+\delta]$. Indeed, since $\cos\beta\le\cos\gamma$, it suffices to show that 
$$
\cos\gamma\le\frac{3\cos(\gamma-\delta)+\cos(\gamma+\delta)}{4}=\cos\gamma\cos\delta+\frac{1}{2}\sin\gamma\sin\delta.
$$
Rewriting this as $\cos\gamma(1-\cos\delta)\le\frac{1}{2}\sin\gamma\sin\delta$ and using the identity $1-\cos\delta$ $=\frac{1-\cos^2\delta}{1+\cos\delta}$ $=\frac{\sin^2\delta}{1+\cos\delta}$, we see that we need to prove that
$\frac{\cos\gamma}{1+\cos\delta}\sin^2\delta\le\frac{1}{2}\sin\gamma\sin\delta$. However, since $0\le\delta\le\gamma\le\frac{\pi}{2}$, we have $\cos\gamma\le\cos\delta\le 1$ and $\sin\gamma\ge\sin\delta$, so the left hand side is at most $\frac{1}{2}\sin^2\delta$ and the right hand side is at least that.

Applying this property to the interval $[\alpha-t_0,\alpha]$, i.e.,  with $\gamma=\alpha-\frac{t_0}{2}$, $\delta=\frac{t_0}{2}$, we conclude that
$$
h_K(e'(t))\le \frac{3}{4}\Bigl(R+\frac{1}{7}|\omega(e)|\Bigr)+\frac{1}{4}(R+\omega(e))=R-\frac{1}{7}|\omega(e)|
$$
for every $t\in [0,\frac{t_0}{2}]\supset [\frac{\vartheta}{5},\frac{2\vartheta}{5}]$ and the conclusion of the lemma follows again.
\ep


\bc\label{FedjaUra}
Let $K$ be a convex body in ${\mathbb R^n}$ and let $R>0$.  Let $h_K=R+\omega$ be the support function of $K$ and let $e\in S^{n-1}$ be a unit vector. Assume that $h_K(e)\le R\cos\vartheta$ for some $\vartheta\in (0,\frac{\pi}{2})$. Then
\[
|\omega(e)|\le C\,\frac{1}{\sigma(S_{\vartheta}(e))}\int_{S_{\vartheta}(e)}|\omega(e')|d\sigma(e').
\]
\ec
\bp
We will use the parametrization $e'=e'(t,v)\in S^{n-1}$ where $t$ is the angle between $e$ and $e'$, and $v\in S^{n-1}\cap e^\perp$ is such that  $e'=e\cos t  + v\sin t $. 

Note that $d\sigma_{n-1}(e')=c_n(\sin t)^{n-2}dtd\sigma_{n-2}(v)$. It follows from the lemma applied to the projection of $K$ to the plane spanned by $e$ and $v$ that
\[
|\omega(e)|\le\frac{35}{\vartheta}\int_{\frac{\vartheta}{5}}^{\vartheta}|\omega(e'(t,v))|dt\le\frac{35}{\vartheta\left(\sin(\frac{\vartheta}{5})\right)^{n-2}}
\int_{0}^{\vartheta}|\omega(e'(t,v))|(\sin t)^{n-2}dt.
\]
Integrating this inequality with respect to $v$ and observing that $\sigma(S_{\vartheta}(e))\asymp\vartheta^{n-1}\asymp\vartheta(\sin(\frac{\vartheta}{5}))^{n-2}$, we get the statement of the corollary.
\ep

\bigskip

\bl\label{rain1}
Assume that a symmetric convex body $K$ is very close to the unit ball and $l\in{\mathbb N}$. Let $h=h_K$ and $\rho=\rho_K$ be the support and the radial functions of $K$ respectively. Trivially, $\rho\le h$. Let $h=\sum\limits_{m=0}^{\infty}h_m$ be the decomposition of $h$ into spherical harmonics (since $h$ is even, only  $h_m$ with even $m$ are not  identically $0$). Put $\eta=\sum\limits_{m=1}^lh_m$, $\nu=\sum\limits_{m=l+1}^{\infty}h_m$. We claim that for every $\varepsilon, l>0$ there exists $\delta_0=\delta_0(\varepsilon, l)$ such that whenever $\|h-1\|_{\infty}\le\delta_0$, the inequality
\[
h-\rho\le \varepsilon\|\eta\|_{L^2}+CM\nu
\]
holds, where $C$ is an absolute constant and $M$ is the spherical Hardy-Littlewood maximal function.
\el
\bp
We have 
\[
\rho(e)=\inf\limits_{\{e'\in S^{n-1},\langle e,e'\rangle>0  \}}\frac{h(e')}{\langle e,e'\rangle}.
\]
Note that the admissible range of $e'$ can be further restricted to $|e-e'|<\delta$ with arbitrarily small $\delta>0$, provided  that $\delta_0$ is chosen small enough. Indeed, since $h(e')\ge \frac{1-\delta_0}{1+\delta_0}h(e)$, $e'$ can compete with $e$ only if $\langle e,e'\rangle\ge \frac{1-\delta_0}{1+\delta_0}$, so 
\[
|e-e'|^2=2(1-\langle e,e'\rangle)\le \frac{4\delta_0}{1+\delta_0}<\delta^2
\]
if $\delta_0>0$ is chosen appropriately. Now observe also that all norms on the finite-dimensional space of polynomials of degree not exceeding $l$ on the unit sphere are equivalent, and that any semi-norm is dominated by any norm, whence
\[
\|\eta\|_{C(S^{n-1})} \le C(l)\|\eta\|_{L^2(S^{n-1})} \qquad\textrm{and}\qquad \|\nabla\eta\|_{C(S^{n-1})}\le C(l)\|\eta\|_{L^2(S^{n-1})}.
\]
In particular, if $|e-e''|<2\delta$, we get 
\[
|\eta(e)-\eta(e'')|\le 4 \|\nabla\eta\|_{C(S^{n-1})}\delta\le 4C(l) \delta \|\eta\|_{L^2(S^{n-1})}.
\]

Let us now assume that $e'\in S^{n-1}$, with $|e-e'|<\delta$, is a competitor, so $\frac{h(e')}{\langle e, e'\rangle}\le h(e)$. Then, if $\vartheta$ is the angle between $e$ and $e'$, we have $h(e')\le h(e)\cos\vartheta$, so we can apply Corollary \ref{FedjaUra} to the vector $e'$ with $R=h(e)$ and conclude that 
\[
h(e)-\frac{h(e')}{\langle e, e'\rangle}\le h(e)-h(e')
\le \frac{C}{\sigma(S_{\vartheta}(e'))}\int\limits_{S_{\vartheta}(e')}|h(e)-h(e'')|d\sigma(e'')
\]
\[
\le \frac{C'}{\sigma(S_{2\vartheta}(e))}\int\limits_{S_{2\vartheta}(e)}|h(e)-h(e'')|d\sigma(e'').
\]
However, 
\[
|h(e)-h(e'')|\le |\eta(e)-\eta(e'')|+|\nu(e)|+|\nu(e'')|,
\]
and
\[
|\eta(e)-\eta(e'')|\le 4C(l)\delta\|\eta\|_{L^2(S^{n-1})},
\]
while 
\[
|\nu(e)|\le M\nu(e)\qquad \textrm{and}\qquad \frac{1}{\sigma(S_{2\vartheta}(e))}\int\limits_{S_{2\theta}(e)}|\nu(e'')|d\sigma(e'')\le M\nu(e),
\]
so the desired statement follows if we choose $\delta>0$ so that $4C'C(l)\delta<\varepsilon$.
\ep

\section{Contraction}

Let ${\mathfrak M}$ be a bounded linear operator on $L^2=L^2(S^{n-1})$ such that ${\mathfrak M}$ is proportional to the identity on every space ${\mathcal H}_m$ of spherical harmonics of degree $m$, i.e.,  for some $\mu_m\in{\mathbb R}$,
\[
{\mathfrak M}f=\sum\limits_{m\ge 0}\mu_mf_m\qquad \textrm{where}\qquad f=\sum\limits_{m\ge 0}f_m
\qquad \textrm{and}\qquad f_m \in {\mathcal H}_m
\]
is the spherical harmonic decomposition of $f\in L^2(S^{n-1})$. We say that ${\mathfrak M}$ is a strong contraction if 
\[
\max\limits_{m\ge 0}|\mu_m|<1,\quad \textrm{and}\qquad \lim\limits_{m\to\infty}\mu_m=0.
\]

\bl\label{rainA}
Assume that ${\mathfrak M}$ as above is a strong contraction. Then, there exists $\delta\in (0,1)$ such that for any symmetric convex body $K$ and any $c\in (1-\delta, 1+\delta)$, the conditions 
\[
1-\delta\le \rho_K\le 1+\delta,\qquad \|(h_K-h_{\scriptscriptstyle{0}})-{c\,\mathfrak M}(\rho_K-r_{\scriptscriptstyle{0}})\|_{L^2}\le \delta \|\rho_K-r_{\scriptscriptstyle{0}}\|_{L^2},
\]
imply $h_K=\rho_K=const$. 
Here, $h_{\scriptscriptstyle{0}}$ and $r_{\scriptscriptstyle{0}}$ are the constant terms of the spherical harmonic decomposition of $h_K$ and $\rho_K$, respectively.
\el
\bp
Fix a large $l$ and consider the decompositions
\[
h_K=h_{\scriptscriptstyle{0}}+\eta+\nu \qquad\textrm{and}\qquad \rho_K=r_{\scriptscriptstyle{0}}+\varphi+\psi,
\]
where $h_{\scriptscriptstyle{0}}$, $r_{\scriptscriptstyle{0}}$ are the constant terms, $\eta$ and $\varphi$ are the parts corresponding to the harmonics of degrees $1$ to $l$  and $\nu$ and $\psi$ are the parts corresponding to the harmonics of degrees greater than $l$.

Fix $\varepsilon>0$. Since the projection to any sum of spaces of spherical harmonics in $L^2$ has norm $1$, we have
\begin{equation}\label{FedjaHot}
\|\eta-c\,{\mathfrak M}\varphi\|_{L^2}\le  \| h_K-h_{\scriptscriptstyle{0}}-c\,{\mathfrak M}(\rho_K-r_{\scriptscriptstyle{0}})\|_{L^2}\le 
\end{equation}
\[
\delta \|\rho_K-r_{\scriptscriptstyle{0}}\|_{L^2} \le \delta(\|\varphi\|_{L^2}+\|\psi\|_{L^2}).
\]
Similarly, 
\begin{equation}\label{Fedjafrost}
\|\nu-c\,{\mathfrak M}\psi\|_{L^2}\le \delta(\|\varphi\|_{L^2}+\|\psi\|_{L^2}).
\end{equation}
From \eqref{Fedjafrost} we obtain that
\[
\|\nu\|_{L^2}\le c\,\|{\mathfrak M}\psi\|_{L^2}+\delta(\|\varphi\|_{L^2}+\|\psi\|_{L^2})\le
\]
\begin{equation}\label{ice}
(1+\delta)(\max\limits_{m> l}|\mu_m|)\|\psi\|_{L^2}+\delta(\|\varphi\|_{L^2}+\|\psi\|_{L^2})\le \varepsilon(\|\varphi\|_{L^2}+\|\psi\|_{L^2})
\end{equation}
if $l$ is large enough (recall that $\lim\limits_{m\to\infty}\mu_m=0$) and $\delta$ is small enough. The same computation for $\eta$, using in this case that   $\max\limits_{m\ge 0}|\mu_m|<1$, yields
\begin{equation}\label{tuz}
   \|\eta\|_{L^2}\le (1+2\delta)(\|\varphi\|_{L^2}+\|\psi\|_{L^2}).
\end{equation}

On the other hand, by Lemma \ref{rain1} and the boundedness of the maximal function in $L^2$, we have 
\[
\|h_K-\rho_K\|_{L^2}\le\varepsilon\|\eta\|_{L^2}+C\|\nu\|_{L^2}, 
\]
which implies 
\begin{equation}\label{do1}
\|\eta-\varphi\|_{L^2}\le\varepsilon\|\eta\|_{L^2}+C\|\nu\|_{L^2}  \quad\textrm{and}\quad \|\nu-\psi\|_{L^2}\le\varepsilon\|\eta\|_{L^2}+C\|\nu\|_{L^2}.
\end{equation}

Combining \eqref{FedjaHot}, \eqref{Fedjafrost}, \eqref{ice}, \eqref{tuz}, and \eqref{do1},  we obtain
\[
\|\varphi-c\,{\mathfrak M}\varphi\|_{L^2}+\|\psi-c\,{\mathfrak M}\psi\|_{L^2} 
\]
\[
\le \|\varphi-\eta\|_{L^2}+ \|\eta -c\,{\mathfrak M}\varphi\|_{L^2} + \|\psi-\nu\|_{L^2} +\|\nu-c\,{\mathfrak M}\psi\|_{L^2} 
\]
\[
\le C(\delta+\varepsilon)(\|\varphi\|_{L^2}+\|\psi\|_{L^2}).
\]
On the other hand, for any function $\chi\in L^2(S^{n-1})$, we have
$$
\|\chi-c\,{\mathfrak M}\chi\|_{L^2}\ge (1-(1+\delta)\max\limits_{m\ge 0}|\mu_m|)\|\chi\|_{L^2},
$$
so we can conclude that 
$\varphi=0$, $\psi=0$ if $\,C(\delta+\varepsilon)<1-(1+\delta)\max\limits_{m\ge 0}|\mu_m|$.
\ep
\br\label{Fedjakrut}
Note that $(h_K-h_{\scriptscriptstyle{0}})-{c\,\mathfrak M}(\rho_K-r_{\scriptscriptstyle{0}})$ is orthogonal to constants and, therefore, its $L^2$-norm does not exceed $\|(h_K-\lambda)-{c\,\mathfrak M}(\rho_K-r_{\scriptscriptstyle{0}})\|_{L^2}$ for any $\lambda \in {\mathbb R}$. Thus, to verify the conditions of the lemma it suffices to check that 
$$
\|(h_K-\lambda)-{c\,\mathfrak M}(\rho_K-r_{\scriptscriptstyle{0}})\|_{L^2}\le \delta \|\rho_K-r_{\scriptscriptstyle{0}}\|_{L^2}
$$
with any $\lambda\in{\mathbb R}$ of our choice.
\er

\section{Properties of the function $(\mathcal{R}[\rho_K^{\alpha}])^{\beta}$ when $\rho_K$ is close to $1$}\label{Fedjaprizemlilsya}

Let $K$ be a symmetric convex body in the isotropic position such that $1-\delta\le\rho_K\le1+\delta$ for some small $\delta>0$. 
Let $\alpha, \beta\in{\mathbb R}$. We want to derive several useful properties of the function $(\mathcal{R}[\rho_K^{\alpha}])^{\beta}$.

The first observation is that $\rho_K$ is Lipschitz with Lipschitz constant $5\sqrt{\delta}$. Indeed, let $x,y\in S^{n-1}$. If $|x-y| \geq \frac{\sqrt{\delta}}{2}$,  then we have 
\[
   |\rho_K(x)-\rho_K(y)| \leq 2 \delta \leq 4 \sqrt{\delta} |x-y|, 
\]
so we may assume that $0<|x-y|<\frac{\sqrt{\delta}}{2}$. Without loss of generality, $\rho_K(x) \geq \rho_K(y)$. Let us denote $X=\rho_K(x)x$, $Y=\rho_K(y)y$, where $X,Y\in \partial K$. By the convexity of $K$, every point on the line $Y-t (X-Y)$ with $t \geq 0$ lies outside $K$ and, therefore, outside $(1-\delta)B_2^n$ as well. Hence,
\[
  (1-\delta)^2 \leq |Y-t (X-Y)|^2 =|Y|^2 -2t\langle X-Y,Y \rangle +t^2|X-Y|^2.
\]
Since $|Y|^2 \leq (1+\delta)^2$, we conclude that, for all $t \geq 0$, 
\begin{equation}\label{lip1}
 2t \langle X-Y,Y \rangle -t^2|X-Y|^2 \leq 4 \delta.
\end{equation}
From \eqref{lip1} it follows that 
\begin{equation}\label{lip2}
    \langle X-Y,Y \rangle \leq 2 \sqrt{\delta} |X-Y|.
\end{equation}
Indeed, if $ \langle X-Y,Y \rangle\le0$, the inequality is obvious. Otherwise, we can plug  $t=\frac{ \langle X-Y,Y \rangle}{|X-Y|^2}$ into \eqref{lip1}, obtaining $\frac{ \langle X-Y,Y \rangle^2}{|X-Y|^2} \leq 4\delta$, which is equivalent to \eqref{lip2}. 
Now, equation \eqref{lip2} can be rewritten as
\[
    |X| |Y| \langle x,y \rangle -|Y|^2 \leq 2 \sqrt{\delta} |X-Y|,
\]
or, equivalently, 
\[
    |Y| \left( |X|- |Y| \right)   \leq 2 \sqrt{\delta} |X-Y| + |X| |Y| (1-\langle x,y \rangle). 
\]
Observe that $1-\langle x,y \rangle=\frac{1}{2}|x-y|^2 \leq \frac{\sqrt{\delta}}{4} |x-y|$, while $|X-Y| \leq |X|-|Y| +|Y||x-y|$. Hence, 
\[
   |X|-|Y| \leq  \frac{2 \sqrt{\delta}}{|Y|} \left( |X| -|Y| \right)+\left( 2+ \frac{|X|}{4} \right) \sqrt{\delta} |x-y|.
\]
Now, if $\delta \in (0,1/25)$, 
\[
\frac{2\sqrt{\delta}}{|Y|} \leq \frac{2\sqrt{\delta}}{1-\delta} \leq \frac{2/5}{1-1/25} <\frac{1}{2},
\]
and we conclude that 
\[
    |X|-|Y| \leq 2 \left(  2 + \frac{|X|}{4} \right) \sqrt{\delta} |x-y| \leq \left( 4 +\frac{1+\delta}{2} \right) \sqrt{\delta} |x-y| \leq 5 \sqrt{\delta} |x-y|,
\]
as required.


Since the mapping $t\mapsto t^p$ is Lipschitz on any compact subset of $(0,+\infty)$ and the Radon transform does not increase the Lipschitz constant of the function, we immediately conclude that $(\mathcal{R}[\rho_K^\alpha])^\beta$ has Lipschitz constant at most $C_{\alpha,\beta}\sqrt{\delta}$.

Let now $r_{\scriptscriptstyle{0}}=
\int\limits_{S^{n-1}}\rho_Kd\sigma$ be the mean value of $\rho_K$ on the unit sphere. Clearly, $|r_{\scriptscriptstyle{0}}-1|\le\delta$, so $|\rho_K-r_{\scriptscriptstyle{0}}|\le 2\delta$ and, thereby, $\mathcal{R}|\rho_K-r_{\scriptscriptstyle{0}}|\le 2\delta$ as well. Now, using the fact that $t\mapsto t^p$ is $C^2$ on any compact 
subset of $(0,+\infty)$ and linearizing, we successively derive that
\[
|\rho_K^\alpha-(r_{\scriptscriptstyle{0}}^\alpha+\alpha r_{\scriptscriptstyle{0}}^{\alpha-1}(\rho_K-r_{\scriptscriptstyle{0}}))|\le C_{\alpha}\delta|\rho_K-r_{\scriptscriptstyle{0}}|,
\]
\[
|\mathcal{R}[\rho_K^\alpha]-(r_{\scriptscriptstyle{0}}^\alpha+\alpha r_{\scriptscriptstyle{0}}^{\alpha-1}\mathcal{R}(\rho_K-r_{\scriptscriptstyle{0}}))|\le C_{\alpha}\,\delta \mathcal{R}|\rho_K-r_{\scriptscriptstyle{0}}|,
\]
\[
|(\mathcal{R}[\rho_K^\alpha])^{\beta}-(r_{\scriptscriptstyle{0}}^\alpha+\alpha r_{\scriptscriptstyle{0}}^{\alpha-1}\mathcal{R}(\rho_K-r_{\scriptscriptstyle{0}}))^{\beta}|\le C_{\alpha, \beta}\,\delta \mathcal{R}|\rho_K-r_{\scriptscriptstyle{0}}|,
\]
\[
|(\mathcal{R}[\rho_K^\alpha])^\beta-(r_{\scriptscriptstyle{0}}^{\alpha\beta}+\alpha\beta r_{\scriptscriptstyle{0}}^{\alpha\beta-1}\mathcal{R}(\rho_K-r_{\scriptscriptstyle{0}}))|\le C_{\alpha,\beta}\,\delta \mathcal{R}|\rho_K-r_{\scriptscriptstyle{0}}|.
\]
Thus, $(\mathcal{R}[\rho_K^\alpha])^\beta=r_{\scriptscriptstyle{0}}^{\alpha\beta}+\gamma$,
where
\[
|\gamma-\alpha\beta r_{\scriptscriptstyle{0}}^{\alpha\beta-1}\mathcal{R}(\rho_K-r_{\scriptscriptstyle{0}})|\le C_{\alpha,\beta}\,\delta \mathcal{R}|\rho_K-r_{\scriptscriptstyle{0}}|.
\]
In particular,  we conclude that $|\gamma|\le C_{\alpha,\beta}\,\delta$, which, together with the above observation about the Lipschitz constant, implies that 
$$
\|\gamma\|_{C^{\frac{1}{2}}}=\max\limits_{S^{n-1}}|\gamma|+\sup\limits_{x,y\in S^{n-1}, x\neq y  }\,\frac{|\gamma(x)-\gamma(y)|}{|x-y|^{\frac{1}{2}}}\le C_{\alpha,\beta}\,\sqrt{\delta}.
$$

Now let $\rho_K-r_{\scriptscriptstyle{0}}=Y_2+Y_4+\dots$ be the spherical harmonic decomposition of 
$\rho_K-r_{\scriptscriptstyle{0}}$. 
It follows from the definition of the isotropic position that
\[
0=\int\limits_Kp(x)dx=c_n\int\limits_{S^{n-1}}\rho_K^{n+2}(x)p(x)d\sigma(x)
\]
for all quadratic polynomials $p(x)=\sum\limits_{i,j}a_{ij}x_ix_j$ with $\sum\limits_{i=1}^na_{ii}=0$. In other words, $\rho_K^{n+2}$ has no second order term in its spherical harmonic decomposition.

On the other hand,
\[
|\rho_K^{n+2}-(
r_{\scriptscriptstyle{0}}^{n+2}+(n+2)r_{\scriptscriptstyle{0}}^{n+1}(\rho_K-r_{\scriptscriptstyle{0}}))|\le C\delta |\rho_K-r_{\scriptscriptstyle{0}}|.
\]
Taking the second order component in the spherical harmonic decomposition of the expression under the absolute value sign on the left hand side, we get
\[
(n+2)r_{\scriptscriptstyle{0}}^{n+1}\|Y_2\|_{L^2(S^{n-1})}\le C\delta \|\rho_K-r_{\scriptscriptstyle{0}}\|_{L^2(S^{n-1})},
\]
so
\[
\|Y_2\|_{L^2(S^{n-1})}\le C'\delta \|\rho_K-r_{\scriptscriptstyle{0}}\|_{L^2(S^{n-1})}.
\]

\section{A  solution to the fifth  Busemann-Petty problem in a small neighborhood of the Euclidean ball}

Recall that for the fifth Busemann-Petty problem we have the equation  $h_K=(\mathcal{R}[\rho_K^{n-1}])^{-1}$. By the results of the previous section, the right hand side can be written as 
\[
r_{\scriptscriptstyle{0}}^{-n+1}-(n-1)
r_{\scriptscriptstyle{0}}^{-n}\mathcal{R}(\rho_K-r_{\scriptscriptstyle{0}})+\gamma',
\]
where $|\gamma'|\le C\delta \mathcal{R}|\rho_K-r_{\scriptscriptstyle{0}}|$, so $\|\gamma'\|_{L^2(S^{n-1})}\le C\delta \|\rho_K-r_{\scriptscriptstyle{0}}\|_{L^2(S^{n-1})}$.

Let ${\mathfrak M}$ be the linear operator that maps every $m$-th order spherical harmonic $Z_m$ to
\[
-(n-1)\mathcal{R}Z_m=-(n-1)(-1)^{\frac{m}{2}}\frac{1\cdot3\cdot\ldots\cdot(m-1)}{(n-1)(n+1)\cdot\ldots\cdot (n+m-3)}Z_m
\]
for even $m\ge 4$ and to $0$ for other $m$. Then ${\mathfrak M}$ is a strong contraction and
\[
\|(h_K-r_{\scriptscriptstyle{0}}^{-n+1})-r_{\scriptscriptstyle{0}}^{-n}{\mathfrak M}(\rho_K-r_{\scriptscriptstyle{0}})\|_{L^2(S^{n-1})}\le
\]
\[
r_{\scriptscriptstyle{0}}^{-n}\|Y_2\|_{L^2(S^{n-1})}+\|\gamma'\|_{L^2(S^{n-1})}\le C\delta \|\rho_K-r_{\scriptscriptstyle{0}}\|_{L^2(S^{n-1})},
\]
so  Lemma \ref{rainA} and Remark \ref{Fedjakrut} yield $h_K=\rho_K=const$, i.e.,  $K$ is a  ball, provided that $\delta$ is small enough.

\section{A  solution to the eighth  Busemann-Petty problem in a small neighborhood of the Euclidean ball}

We now turn to the equation 
$Ah_K=(\mathcal{R}[\rho_K^{n-1}])^{n+1}$ (see Section \ref{Fedjaletaet}). Below we will use several standard results about $A$ and the Laplace operator which, for completeness, are proven in the Appendices. 

By the results of Section \ref{Fedjaprizemlilsya}, $(\mathcal{R}[\rho_K^{n-1}])^{n+1}$ can be rewritten as $r_{\scriptscriptstyle{0}}^{(n-1)(n+1)}+\gamma$, where 
$\|\gamma\|_{C^{\frac{1}{2}}}\le C\sqrt{\delta}$ and
\[
\gamma=(n-1)(n+1)r_{\scriptscriptstyle{0}}^{(n-1)(n+1)-1}\mathcal{R}(\rho_K-r_{\scriptscriptstyle{0}})+\gamma',
\]
\[
\|\gamma'\|_{L^2(S^{n-1})}\le C\delta\|\rho_K-r_{\scriptscriptstyle{0}}\|_{L^2(S^{n-1})}.
\]
Then
\[
A\frac{h_K}{r_{\scriptscriptstyle{0}}^{n+1}}=1+r_{\scriptscriptstyle{0}}^{-(n-1)(n+1)}\gamma
\]
and, provided that $\delta>0$ is small enough, we can apply Lemma \ref{Fl4} (see Appendix II) 
and the uniqueness theorem (see \cite{Sch}, Theorem 8.1.1)
to obtain 
\[
\frac{h_K}{r_{\scriptscriptstyle{0}}^{n+1}}=1+\varphi'+\varphi'',
\]
where 
\begin{equation}\label{dob113}
	\widetilde{\Delta}\varphi'=r_{\scriptscriptstyle{0}}^{-(n-1)(n+1)}\gamma
\end{equation}
(see Appendix II for the definition of $\widetilde{\Delta}$) and 
\[
 \|\varphi''\|_{L^2(S^{n-1})}\le \varepsilon r_{\scriptscriptstyle{0}}^{-(n-1)(n+1)} \|\gamma\|_{L^2(S^{n-1})}\le C\varepsilon \|\rho_K-r_{\scriptscriptstyle{0}}\|_{L^2(S^{n-1})}
\]
with as small $\varepsilon>0$ as we want.

Furthermore, the solution of  equation \eqref{dob113} splits into $\varphi'_1+\varphi_2'$ where $\varphi_1'$ solves 
\[
	\widetilde{\Delta}\varphi_1'=(n-1)(n+1)r_{\scriptscriptstyle{0}}^{-1}{\mathcal R}(\rho_K-r_{\scriptscriptstyle{0}}),
\]
and
$\varphi_2'$ solves $\widetilde{\Delta}\varphi_2'=r_{\scriptscriptstyle{0}}^{-(n-1)(n+1)}\gamma'$.

The norm of $\varphi_2'$ can be estimated immediately:
\[
\|\varphi_2'\|_{L^2(S^{n-1})}\le C \|\gamma'\|_{L^2(S^{n-1})}\le C\delta \|\rho_K-r_{\scriptscriptstyle{0}}\|_{L^2(S^{n-1})}.
\]
As to $\varphi_1'$, it is equal to (see the end of Appendix I)
\[
r_{\scriptscriptstyle{0}}^{-1}\sum\limits_{{\substack{m\ge 2\\ m\,\text{even}}}}\mu_mY_m,
\]
where
$\rho_K=r_{\scriptscriptstyle{0}}+\sum\limits_{{\substack{m\ge 2\\ m\,\text{even}}}}Y_m\,$
is the spherical harmonic decomposition of $\rho_K$ and
\[
\mu_m=\frac{(n-1)(n+1)}{(1-m)(m+n-1)}(-1)^{\frac{m}{2}}\frac{1\cdot3\cdot\ldots\cdot(m-1)}{(n-1)(n+1)\cdot\ldots\cdot (n+m-3)},
\]
so $\mu_2=1$ and $|\mu_m|<1$ for $m\ge 4$, $\mu_m\to 0$ as $m\to\infty$. Since
\[
\|Y_2\|_{L^2(S^{n-1})}\le  C\delta \|\rho_K-r_{\scriptscriptstyle{0}}\|_{L^2(S^{n-1})},
\]
we conclude that
\[
\|\varphi_1'-r_{\scriptscriptstyle{0}}^{-1}{\mathfrak M}(\rho_K-r_{\scriptscriptstyle{0}})\|_{L^2(S^{n-1})}\le C\delta\|\rho_K-r_{\scriptscriptstyle{0}}\|_{L^2(S^{n-1})},
\]
with the strong contraction ${\mathfrak M}$ given by
$Z_m\mapsto\mu_mZ_m$, $m$ even, $m\ge 4$; $Z_m\mapsto 0$ for all other $m$.

Putting all these estimates together, we conclude that
\[
\|(h_K-r_{\scriptscriptstyle{0}}^{n+1})-r_{\scriptscriptstyle{0}}^{n}{\mathfrak M}(\rho_K-r_{\scriptscriptstyle{0}})\|_{L^2(S^{n-1})}\le \varepsilon\|\rho_K-r_{\scriptscriptstyle{0}}\|_{L^2(S^{n-1})},
\]
with as small $\varepsilon>0$ as we want, provided that $\delta>0$ is small enough. Now we can apply Lemma \ref{rainA} and Remark \ref{Fedjakrut} again to conclude that $h_K=\rho_K=const$, so $K$ is a ball.

\section{Appendix I. Solving the Laplace equation}

Below we shall use the following notation. For a function $f:S^{n-1}\to{\mathbb R}$ and $\alpha\in (0,1)$, we shall denote
\[
\|f\|_{C^{\alpha}}=\|f\|_{C^{\alpha}(S^{n-1})}=\max\limits_{S^{n-1}}|f|+\sup\limits_{x,y\in S^{n-1}, \,x\neq y}  \frac{|f(x)-f(y)|}{|x-y|^{\alpha}},
\]
\[
\|f\|_{C^{2+\alpha}}=\|f\|_{C^{2+\alpha}(S^{n-1})}=
\]
\[
\max\limits_{S^{n-1}}|f|+\max\limits_{x\in S^{n-1},\,i=1,\dots,n}|F_{x_i}(x)|+
\max\limits_{i,j =1,\dots,n}\|F_{x_ix_j}\|_{C^{\alpha}(S^{n-1})},
\]
where
$F(x)=|x|f(\frac{x}{|x|})$ is the $1$-homogeneous extension of $f$ to ${\mathbb R^n}\setminus\{0\}$ (we assume that it is at least $C^2$ in ${\mathbb R^n}\setminus\{0\}$).

Let $g:S^{n-1}\to{\mathbb R}$ be an even $C^\alpha$ function on the unit sphere $S^{n-1}$ with some $\alpha\in (0,1)$. Let ${\mathcal G}$ be the $(-1)$-homogeneous extension of $g$ to ${\mathbb R^n}\setminus \{0\}$, i.e., ${\mathcal G}(x)=|x|^{-1}g(\frac{x}{|x|})$ for $x\neq 0$.  We will show that there exists a unique $1$-homogeneous even function $F:{\mathbb R^n}\to{\mathbb R}$ of class $L^1_{loc}$ such that  $\Delta F={\mathcal G}$ in ${\mathbb R^n}$ in the sense of generalized functions. Moreover, $F\in C^{2+\alpha}(S^{n-1})$ and for all $i,j=1,\dots,n$, we have
\begin{equation}\label{dob1}
\|F_{x_ix_j}\|_{L^2(S^{n-1})}\le C \|g\|_{L^2(S^{n-1})},\qquad\|F\|_{C^{2+\alpha}(S^{n-1})}\le C \|g\|_{C^{\alpha}(S^{n-1})},
\end{equation}
with some $C=C(n,\alpha)>0$.

\subsection{Uniqueness} If we have two even $1$-homogeneous functions $F_1$, $F_2$ such that $\Delta F_1=\Delta F_2={\mathcal G}$  in ${\mathbb R^n}$, then $F_1-F_2$ is an even $1$-homogeneous harmonic function, but the only such function is $0$.

\subsection{Existence} Now we will show that the function $F$ defined by 
\[
F(x)=c_n\int\limits_{\mathbb R^n}\Big[\frac{1}{|x-y|^{n-2}}-\frac{1}{|y|^{n-2} }  \Big]{\mathcal G}(y)dy
\]
is a well-defined $1$-homogeneous function on ${\mathbb R^n}$ satisfying $\Delta  F={\mathcal G}$ and estimates \eqref{dob1}. Here, $c_n$ is chosen so that $\Delta \frac{c_n}{|x|^{n-2}}=\delta_0$ (the Dirac delta measure) in the sense of generalized functions, and the integral is understood as $\lim\limits_{R\to\infty}\int\limits_{B(0,R)}$.

 In order to show the convergence of the integral, we note that
\[
\frac{1}{|x-y|^{n-2}}-\frac{1}{|y|^{n-2}}=(n-2)\frac{\langle x,y\rangle}{|y|^{n}}+O\Big(\frac{1}{|y|^n}\Big)
\]
as $y\to\infty$ uniformly on compact sets in $x$.

Since ${\mathcal G}$ is even, the integral of $\frac{\langle x,y\rangle}{|y|^{n}}{\mathcal G}(y)$ over each sphere centered at the origin vanishes. Since ${\mathcal G}$ is $(-1)$-homogeneous, we have $\frac{1}{|y|^{n}}{\mathcal G}(y)=O(\frac{\|g\|_{C}}{|y|^{n+1}})$ as $y\to \infty$, which is integrable at $\infty$.

The singularities at  $x\in S^{n-1}$ and $0$ are of degrees $-(n-2)$ and $-(n-1)$   respectively, so the local integrability there presents no problem either, and we get the estimate
\[
\|F\|_{C(S^{n-1})}\le C\|g\|_{C(S^{n-1})}.
\]

The change of variable $y\mapsto -y$ and the identity ${\mathcal G}(y)={\mathcal G}(-y)$ imply that $F$ is even. 

To show  the $1$-homogeneity of $F$, take $t>0$ and apply the change of variable $y\mapsto ty$ to write
\[
F(tx)=c_n\lim\limits_{R\to\infty}\int\limits_{B(0,R)}\Big(\frac{1}{|tx-y|^{n-2}}- \frac{1}{|y|^{n-2}}      \Big){\mathcal G}(y)dy=
\]
\[
c_n\lim\limits_{R\to\infty}\int\limits_{B(0,\frac{R}{t})}\Big(\frac{1}{|tx-ty|^{n-2}}- \frac{1}{|ty|^{n-2}}      \Big){\mathcal G}(ty)d(ty)=
\]
\[
c_n\,t\,\lim\limits_{R\to\infty}\int\limits_{B(0,\frac{R}{t})}\Big(\frac{1}{|x-y|^{n-2}}- \frac{1}{|y|^{n-2}}      \Big){\mathcal G}(y)dy=tF(x),
\]
(we used that $\frac{1}{|tz|^{n-2}}=\frac{1}{t^{n-2}}\frac{1}{|z|^{n-2}}$, ${\mathcal G}(ty)=t^{-1}{\mathcal G}(y)$, and $d(ty)=t^ndy$).

\begin{figure}[h]
	\centering
	\includegraphics[height=1in]{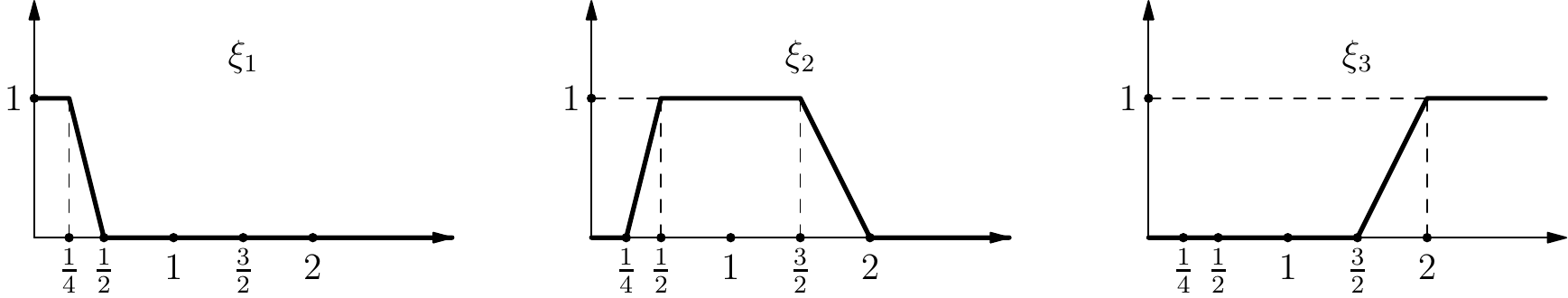} 
	\caption{The functions $\xi_i$, $i=1,2,3$}
	\label{ksi}
\end{figure}

To estimate $\|F\|_{C^{2+\alpha}(S^{n-1})}$, we split the integral defining $F$  into 3 parts. Let 
$\xi_1$, 
$\xi_2$, 
$\xi_3:[0,+\infty)$
$\to [0,1]$ be as on  Figure \ref{ksi}, so $\xi_i$ are Lipschitz with constant 4, and $\xi_1+\xi_2+\xi_3=1$.

Put ${\mathcal G}_i(x)={\mathcal G}(x)\xi_i(|x|)$ and 
\[
F_i(x)=c_n\int\limits_{{\mathbb R^n}}\Big[\frac{1}{|x-y|^{n-2}}- \frac{1}{|y|^{n-2}}   \Big]{\mathcal G}_i(y)dy,
\]
so ${\mathcal G}={\mathcal G}_1+{\mathcal G}_2+{\mathcal G}_3$ and $F=F_1+F_2+F_3$.

Our first observation is that 
 ${\mathcal G}_2(y)$ is an $\alpha$-H\"older,  compactly supported function on ${\mathbb R^n}$, with  $C^\alpha$-norm bounded by $C\|g\|_{C^{\alpha}(S^{n-1})}$.

Indeed, we clearly have $\max\limits_{\mathbb R^n}|{\mathcal G}_2|\le 4\max\limits_{S^{n-1}}|g|$.
On the other hand, 
\[
|{\mathcal G}_2(x)-{\mathcal G}_2(y)|=\Big|\xi_2(|x|)|x|^{-1}g\Big(\frac{x}{|x|}\Big)-\xi_2(|y|)|y|^{-1}g\Big(\frac{y}{|y|}\Big)\Big|.
\]
Since $\widetilde{\xi}_2(t)=\xi_2(t)t^{-1}$ is a compactly supported Lipschitz function on $[0,+\infty)$, it is also $\alpha$-H\"older for any $\alpha\in (0,1)$, i.e., 
\[
|\widetilde{\xi}_2(t)-\widetilde{\xi}_2(s)|\le C|t-s|^{\alpha}\quad \textrm{for all}\quad t,s\ge 0.
\]
Thus, if $x,y\in \overline{B(0,2)}\setminus B(0,\frac{1}{4})$, then
\[
\left|\widetilde{\xi}_2(|x|)g\left(\frac{x}{|x|}\right)-\widetilde{\xi}_2(|y|)g\left(\frac{y}{|y|}\right)\right|\le
\]
\[
|\widetilde{\xi}_2(|x|)-\widetilde{\xi}_2(|y|)|\,\left|g\left(\frac{x}{|x|}\right)\right|+|\widetilde{\xi}_2(|y|)|\, \left|g\left(\frac{x}{|x|}\right)-g\left(\frac{y}{|y|}\right)\right|\le
\]
\[
C\,\Bigl||x|-|y|\Bigr|^{\alpha}\max\limits_{S^{n-1}}|g|+4\Big|\frac{x}{|x|}-\frac{y}{|y|}\Big|^{\alpha}\|g\|_{C^{\alpha}(S^{n-1})}\le
\]
\[
C\,\|g\|_{C^{\alpha}(S^{n-1})}\Big(|x-y|^{\alpha}+\Big|\frac{x}{|x|}-\frac{y}{|y|}\Big|^{\alpha}\Big)\le 
C\,\|g\|_{C^{\alpha}(S^{n-1})}|x-y|^{\alpha},
\]
because the mapping $x\mapsto \frac{x}{|x|}$ is $C^1$ and, thereby, Lipschitz on $\overline{B(0,2)}\setminus B(0,\frac{1}{4})$.

If $x,y \notin\overline{B(0,2)}\setminus B(0,\frac{1}{4})$, then ${\mathcal G}_2(x)={\mathcal G}_2(y)=0$, so the inequality 
\[
|{\mathcal G}_2(x)-{\mathcal G}_2(y)|\le C\,\|g\|_{C^{\alpha}(S^{n-1})}|x-y|^\alpha
\]
holds trivially. Finally, if $x\in \overline{B(0,2)}\setminus B(0,\frac{1}{4})$ but $y\notin \overline{B(0,2)}\setminus B(0,\frac{1}{4})$, then the segment $[x,y]$ intersects the boundary of $\overline{B(0,2)}\setminus B(0,\frac{1}{4})$ at some point $y'$, so
${\mathcal G}_2(y)={\mathcal G}_2(y')=0$ and
\[
|{\mathcal G}_2(x)-{\mathcal G}_2(y)|= |{\mathcal G}_2(x)-{\mathcal G}_2(y')|\le
\]
\[
 C\,\|g\|_{C^{\alpha}(S^{n-1})}|x-y'|^{\alpha}
\le C\,\|g\|_{C^{\alpha}(S^{n-1})}|x-y|^{\alpha}.
\]
The functions ${\mathcal G}_1$ and ${\mathcal G}_3$ are supported on $\overline{B(0,\frac{1}{2})}$ and 
${\mathbb R^n}\setminus B(0,\frac{3}{2})$, respectively, and satisfy the bound 
\[
|{\mathcal G}_1(y)|, \,|{\mathcal G}_3(y)|\le \frac{1}{|y|}\max\limits_{S^{n-1}}|g|.
\]

Now we are ready to estimate  $\|F\|_{C^{2+\alpha}(S^{n-1})}$. Consider $x$ with $\frac{3}{4}\le|x|\le\frac{5}{4}$. Note that $x\mapsto\frac{1}{|x-y|^{n-2}}$ is a $C^3$-function (in $x$) in this domain with uniformly bounded (in $y$) $C^3$-norm as long as $y\in \overline{B(0,\frac{1}{2})}$. Hence, 
$F_1\in C^3(\overline{B(0,\frac{5}{4})}\setminus B(0,\frac{3}{4}))$ and 
\[
\|F_1\|_{C^3(\overline{B(0,\frac{5}{4})}\setminus B(0,\frac{3}{4}))}\le C\|g\|_{L^1(S^{n-1})}
\]
(the constant term $\int\limits_{{\mathbb R^n}} \frac{1}{|y|^{n-2}}{\mathcal G}_1(y)dy$ is also bounded by 
$C\|g\|_{L^1(S^{n-1})}$).

To estimate $F_3$, note that for $|x|\le\frac{5}{4}$ and $|y|\ge \frac{3}{2}$, we have
\[
\Big|\frac{1}{|x-y|^{n-2}}- \frac{1}{|y|^{n-2}}  -(n-2)\frac{\langle x,y\rangle}{|y|^{n}}    \Big|\le \frac{C}{|y|^{n}},
\]
\[
\Big|\frac{\partial}{\partial x_i}\frac{1}{|x-y|^{n-2}} -(n-2)\frac{y_i}{|y|^{n}}    \Big|\le \frac{C}{|y|^{n}},
\]
\[
\Big|\frac{\partial^2}{\partial x_i\partial x_j}\frac{1}{|x-y|^{n-2}}   \Big|\le \frac{C}{|y|^{n}},
\]
and
\[
\Big|\frac{\partial^3}{\partial x_i\partial x_j\partial x_k}\frac{1}{|x-y|^{n-2}}   \Big|\le \frac{C}{|y|^{n+1}}.
\]

Since $y \mapsto \frac{\langle x,y\rangle}{|y|^{n}} $ and $y\mapsto 
\frac{y_i}{|y|^{n}} $ are odd functions, their integrals against the even function ${\mathcal G}_3(y)$ over any sphere centered at the origin are $0$ and, therefore, 
\[
|F_3|,|\nabla F_3|, |\nabla^2 F_3|\le C\int\limits_{{\mathbb R^n}\setminus B(0,\frac{3}{2})}|y|^{-n}|{\mathcal G}_3(y)|dy\le C\|g\|_{L^1(S^{n-1})}
\]
and
\[
 |\nabla^3 F_3|\le C\int\limits_{{\mathbb R^n}\setminus B(0,\frac{3}{2})}|y|^{-n-1}|{\mathcal G}_3(y)|dy\le C\|g\|_{L^1(S^{n-1})},
\]
so
\[
\|F_3\|_{C^3(\overline{B(0,\frac{5}{4})})}\le C\|g\|_{L^1(S^{n-1})}.
\]

It remains to estimate $F_2$.
We clearly have
\[
|F_2(x)|\le C\|g\|_{C(S^{n-1})}\int\limits_{\overline{B(0,2)}\setminus B(0,\frac{1}{4})}\Big|\frac{1}{|x-y|^{n-2}}- \frac{1}{|y|^{n-2}}     \Big|dy\le C\|g\|_{C(S^{n-1})}
\]
and
\[
|\nabla F_2(x)|\le C\|g\|_{C(S^{n-1})}\int\limits_{\overline{B(0,2)}\setminus B(0,\frac{1}{4})}\frac{1}{|x-y|^{n-1}}dy\le C\|g\|_{C(S^{n-1})}.
\]

As for $(F_2)_{x_ix_j}$, these partial derivatives are images
of ${\mathcal G}_2$ under certain Calder\'on-Zygmund singular integral operators (see \cite{GT}, Lemma 4.4 and Theorem 9.9), so, since 
${\mathcal G}_2\in C^\alpha{({\mathbb R^n})}$ and has fixed compact support, we obtain that
\[
\|(F_2)_{x_ix_j}\|_{C^{\alpha}({\mathbb R^n})}\le C\|{\mathcal G}_2\|_{C^{\alpha}({\mathbb R^n})}\le C\|g\|_{C^{\alpha}(S^{n-1})}
\]
and
\[
\|(F_2)_{x_ix_j}\|_{L^2({\mathbb R^n})}\le C\|{\mathcal G}_2\|_{L^2({\mathbb R^n})}\le C\|g\|_{L^2(S^{n-1})}.
\]

The final conclusion is that
\[
\|F\|_{C^{2+\alpha}(S^{n-1})}\le C\|F\|_{C^{2+\alpha}(\overline{B(0,\frac{5}{4})}\setminus B(0,\frac{3}{4}))}\le C\|g\|_{C^{\alpha}(S^{n-1})}
\]
and
\[
\|(F_2)_{x_ix_j}\|_{L^2(S^{n-1})}\le C\|(F_2)_{x_ix_j}\|_{L^2(\overline{B(0,\frac{5}{4})}\setminus B(0,\frac{3}{4}))}\le C\|g\|_{L^2(S^{n-1})}
\]
(we used the $(-1)$-homogeneity of $(F_2)_{x_ix_j}$  here).

The desired equality $\Delta F={\mathcal G}$ follows from the fact that the mapping $x\mapsto \frac{1}{|x-y|^{n-2}}-\frac{1}{|y|^{n-2}}$ is harmonic in $x$ for $|x|\le 1$, $|y|\ge \frac{3}{2}$. This implies that 
$\Delta F_3=0$ in $B(0,1)$, while $F_1+F_2$ differs by a constant from the classical Newton potential of the compactly supported $L^1$ function 
${\mathcal G}_1+{\mathcal G}_2={\mathcal G}$ in $B(0,1)$. Hence,
 $\Delta F={\mathcal G}$ in $B(0,1)$, 
and this identity extends to ${\mathbb R^n}$ by homogeneity.

We shall also need the relation between the spherical harmonic decompositions of $F|_{S^{n-1}}$ and $g$.
To this end, we will start with the following computation. Let $P_m$ be a homogeneous harmonic polynomial of degree $m$, so that $Y_m=P_m|_{S^{n-1}}$ is a spherical  harmonic of degree $m$. The $1$-homogeneous extension of $Y_m$ is $\widetilde{Y}_m(x)=|x|^{1-m}P_m(x)$. Then
\[
\Delta \widetilde{Y}_m(x)=\Delta(|x|^{1-m}) P_m(x)+2\langle \nabla (|x|^{1-m}),\nabla P_m(x)\rangle=
\]
\[
(1-m)(-m-1+n)|x|^{-m-1}P_m(x)+2(1-m)|x|^{-m}\frac{\partial}{\partial r}P_m(x)=
\]
\[
|x|^{-1-m}(1-m)(-m-1+n+2m)P_m(x)=
\]
\[
(1-m)(m+n-1)|x|^{-1-m}P_m(x).
\]
Thus, if $g\in L^2(S^{n-1})$ and $g=\sum\limits_{{\substack{m\ge 0\\ m\,\text{even}}}}Y_m$ on $S^{n-1}$,  the series 
\[
F=\sum\limits_{{\substack{m\ge 0\\ m\,\text{even}}}}\frac{1}{(1-m)(m+n-1)}\widetilde{Y}_m
\]
converges in $L^2_{loc}$ (the series is orthogonal on every ball $B(0,R)$ and 
$\|\widetilde{Y}_m\|_{L^2(B(0,R))}$ $\le C_R\|Y_m\|_{L^2(S^{n-1})})$
and formally solves $\Delta F={\mathcal G}$. To show that it is a true solution, it suffices to observe that we have $\Delta F_{(l)}={\mathcal G}_{(l)}$ for the partial sums $F_{(l)}$ and ${\mathcal G}_{(l)}$ of the corresponding series and $F_{(l)}\to F$ in $L^2_{loc}$, ${\mathcal G}_{(l)}\to{\mathcal G}$ in $L^1_{loc}$ as $l\to\infty$. Thus, $\Delta F={\mathcal G}$ in the sense of generalized functions. If $g\in C^\alpha(S^{n-1})$, then, by the uniqueness part, this solution has to coincide with the explicit solution constructed above, so the spherical harmonic decomposition of $F|_{S^{n-1}}$ is $\sum\limits_{\substack{m\ge 0\\ m\,\text{even}}}\frac{1}{(1-m)(m+n-1)}Y_m$.
In particular, the decomposition implies that 
\[
\|F\|_{L^2(S^{n-1})}\le \|g\|_{L^2(S^{n-1})}.
\]

 \section{Appendix II. Solution of Monge-Ampere equation}
 
 For a function $f:S^{n-1}\to{\mathbb R}$, we denote by $F$ its $1$-homogeneous extension to ${\mathbb R^n}$.
 By $Af$ we will denote the restriction of $\sum\limits_{k=1}^n\det \widehat{F}_k$ to the unit sphere where 
 $\widehat{F}_k$ is the matrix obtained from the Hessian $\widehat{F}=(F_{x_ix_j})^n_{i,j=1}$ by deleting the $k$-th row and the $k$-th column.

 We now turn to the solution of the equation $Af=g$ where $g$ is close to $1$. 
 Note that $A1=1$. Indeed, since $A$ commutes with the rotations of the sphere, we can check this identity at the point $(1,0,\dots,0)$. The $1$-homogeneous extension of $1$ is $|x|$, so the Hessian is $\Big(\frac{\delta_{ij}}{|x|}-\frac{x_ix_j}{|x|^3} \Big)_{i,j=1}^n$, which at the point $(1,0,\dots,0)$ turns into 
 $$
 \begin{bmatrix} 0 & 0  & 0 & \dots & 0  & 0\\  0 & 1 & 0 & \dots &0 &0  \\ 0 & 0 & 1 & \dots &0  &0\\
 \vdots & \vdots & \vdots & \ddots & \vdots & \vdots \\ 0 & 0 & 0 & \dots &1  &0 \\ 0 & 0  & 0 & \dots & 0  & 1\end{bmatrix}.
 $$
 
The rotation invariance also allows us to compute the linear part of $Af$  (meaning the linear terms in $\Phi_{x_ix_j}$) for $f=1+\varphi$, where $\Phi$ is the $1$-homogeneous extension of $\varphi$. Again, computing the Hessian at $(1,0,\dots,0)$, we get
 $$
 \widehat{F}=\begin{bmatrix} \Phi_{x_1x_1} &  \Phi_{x_1x_2}  & \dots &   \Phi_{x_1x_n}\\  \Phi_{x_2x_1} &  1+\Phi_{x_2x_2}  & \dots &   \Phi_{x_2x_n}  \\
 \vdots & \vdots & \ddots & \vdots  \\\Phi_{x_nx_1} &  \Phi_{x_nx_2}  & \dots &   1+\Phi_{x_nx_n} \end{bmatrix},
 $$
 so
\[
 \sum\limits_{i=1}^n \det \widehat{F}_i=\det \widehat{F}_1+\sum\limits_{i=2}^n \det \widehat{F}_i=
\]
\[
 1+\sum\limits_{i=2}^n\Phi_{x_ix_i}+(n-1)\Phi_{x_1x_1}+P(\Phi),
\]
 where $P(\Phi)$ is some linear combination of products of two or more second partial derivatives of $\Phi$. Note now that, since $\Phi$ is $1$-homogeneous, the mapping $t\mapsto \Phi(t,0,\dots,0)$ is linear and, thereby, $\Phi_{x_1x_1}(1,0,\dots,0)=0$. Thus we can just as well write
 $\sum\limits_{i=2}^n\Phi_{x_ix_i}+(n-1)\Phi_{x_1x_1}$
 at $(1,0,\dots,0)$ as $\Delta\Phi(1,0,\dots,0)$. However, $\Delta\Phi$ also commutes with rotations, so we  have the identity
\[
 \sum\limits_{i=1}^n \det \widehat{F}_i=1+\Delta\Phi+P(\Phi)
\]
 in general, though $P(\Phi)$ will now be a sum of products of at least two second partial derivatives of $\Phi$ and some fixed   functions of $x$ that are    smooth
 near the unit sphere.
 
 Using identities of the type
\[
 a_1a_2\dots a_m-b_1b_2\dots b_m=(a_1-b_1)a_2\dots a_m+
\]
\[
 b_1(a_2-b_2)a_3\dots a_m+ \cdots + b_1\dots b_{m-2}(a_{m-1}-b_{m-1})a_m+b_1\dots b_{m-1}(a_m-b_m),
\]
 we see that for any $1$-homogeneous $C^2$-functions $\Psi'$, $\Psi''$ satisfying 
\[
 \max\limits_{i,j} \|\Psi_{x_ix_j}'\|_{C^\alpha(S^{n-1})}\le 1,\quad
 \max\limits_{i,j} \|\Psi_{x_ix_j}''\|_{C^\alpha(S^{n-1})}\le 1,
\]
 we have
 \begin{equation}\label{vc1}
 \|P(\Psi')-P(\Psi'')\|_{L^2(S^{n-1})}\le 
 \end{equation}
\[
C\max\limits_{i,j} \|\Psi_{x_ix_j}'-\Psi_{x_ix_j}''\|_{L^2(S^{n-1})}  \max\limits_{i,j} \Bigl(\|\Psi_{x_ix_j}'\|_{C(S^{n-1})}+\|\Psi_{x_ix_j}''\|_{C(S^{n-1})}\Bigr)
\]
 and 
 \begin{equation}\label{vc2}
 \|P(\Psi')-P(\Psi'')\|_{C^{\alpha}(S^{n-1})}\le 
 \end{equation}
\[
C\max\limits_{i,j} \|\Psi_{x_ix_j}'-\Psi_{x_ix_j}''\|_{C^{\alpha}(S^{n-1})}  \max\limits_{i,j} \Bigl(\|\Psi_{x_ix_j}'\|_{C^\alpha(S^{n-1})}+\|\Psi_{x_ix_j}''\|_{C^\alpha(S^{n-1})}\Bigr).
\]
 This will enable us to solve the equation 
 $Af=g$ with $g=1+\gamma$ by iterations if $\|\gamma\|_{C^{\alpha}(S^{n-1})}$ is small enough.

 By $\widetilde{\Delta} f$ we shall denote the restriction of the Laplacian $\Delta F$ of $F$ to the unit sphere. Note that the Laplacian $\Delta F$ itself is a $(-1)$-homogeneous function on ${\mathbb R^n}\setminus\{0\}$ (assuming again that $F$ is twice continuously differentiable away from the origin).


 \bl\label{Fl4} 
 For every $\varepsilon>0$, there exists $\delta>0$ such that for every even $g=1+\gamma$ with $\|\gamma\|_{C^{\alpha}(S^{n-1})}\le \delta$, there exists $f=1+\varphi$ solving $Af=g$ and such that $\|\varphi\|_{C^{2+\alpha}(S^{n-1})}\le\varepsilon$ and, moreover,  $\varphi=\varphi'+\varphi''$, where
 $\widetilde{\Delta} \varphi'=\gamma$, while $ \|\varphi''\|_{L^2(S^{n-1})}\le \varepsilon \|\gamma\|_{L^2(S^{n-1})}$.
 \el
 \bp
 Define the sequence $\varphi_m$ as follows: $\widetilde{\Delta}\varphi_0=\gamma$, $\widetilde{\Delta}\varphi_1=\gamma-P(\Phi_0)$, 
 $\widetilde{\Delta}\varphi_2=\gamma-P(\Phi_1)$, etc., where as before, $\Phi_m$ is the $1$-homogeneous extension of $\varphi_m$.  Recall that by the results of Appendix I, for every even  function $\chi\in C^{\alpha}(S^{n-1})$, there exists a unique solution $\psi$ of the equation $\widetilde{\Delta}\psi=\chi$ and we have the estimates
 $$
 \|\psi\|_{C^{2+\alpha}(S^{n-1})}\le K\|\chi\|_{C^{\alpha}(S^{n-1})},\quad \max\limits_{i,j}\|\Psi_{x_ix_j}\|_{L^2(S^{n-1})}\le K\|\chi\|_{L^2(S^{n-1})}
 $$
 with some  constant $K>0$. So, all $\varphi_m$ are well-defined.

 Let $\delta>0$ be a very small number. Then,  under the assumption 
 $\|\gamma\|_{C^{\alpha}(S^{n-1})}\le \delta$, we have
 $$
 \|\varphi_0\|_{C^{2+\alpha}(S^{n-1})}\le K\|\gamma\|_{C^{\alpha}(S^{n-1})}\le K\delta.
 $$
  Fix $\kappa\in(0,\frac{1}{K})$. It follows from \eqref{vc2} that as long as
 $$
 \|\psi'\|_{C^{2+\alpha}(S^{n-1})},\,\|\psi''\|_{C^{2+\alpha}(S^{n-1})}\le\frac{\kappa}{2C},
 $$
 we have
 $$
 \|P(\Psi')-P(\Psi'')\|_{C^{\alpha}(S^{n-1})}\le\kappa \|\psi'-\psi''\|_{C^{2+\alpha}(S^{n-1})}.
 $$
 If $\delta$ is small enough, so that $K\delta<\frac{\kappa}{2C}$, we obtain
 $$
  \|P(\Phi_0)\|_{C^{\alpha}(S^{n-1})}=\|P(\Phi_0)-P(0)\|_{C^{\alpha}(S^{n-1})}\le \kappa K\delta.
 $$
 Hence, from
 $\widetilde{\Delta}(\varphi_1-\varphi_0)=-P(\Phi_0)$, we conclude that
 $$
 \|\varphi_1-\varphi_0\|_{C^{2+\alpha}(S^{n-1})}\le K(\kappa K)\delta,
 $$
 so
 $$
 \|\varphi_1\|_{C^{2+\alpha}(S^{n-1})}\le K(1+\kappa K)\delta.
 $$
 If this value is still less that $\frac{\kappa}{2C}$, we can continue and write
 $$
 \|P(\Phi_1)-P(\Phi_0)\|_{C^{\alpha}(S^{n-1})}\le \kappa \|\varphi_1-\varphi_0\|_{C^{2+\alpha}(S^{n-1})}\le (\kappa K)^2\delta.
 $$
 Thus, from $\widetilde{\Delta}(\varphi_2-\varphi_1)=-(P(\Phi_1)-P(\Phi_0))$, we get
 $$
 \|\varphi_2-\varphi_1\|_{C^{2+\alpha}(S^{n-1})}\le K (\kappa K)^2\delta,
 $$
 $$
 \|\varphi_2\|_{C^{2+\alpha}(S^{n-1})}\le K(1+\kappa K+(\kappa K)^2)\delta,
 $$
 and so on. We can continue this chain of estimates as long as
  $$
 K(1+\kappa K+(\kappa K)^2+\dots+(\kappa K)^m)\delta<\frac{\kappa}{2C}, 
 $$
 which is forever  if $\frac{K\delta}{1-\kappa K}<\frac{\kappa}{2C}$. 
 
 The outcome is that 
 $$
  \|\varphi_{m+1}-\varphi_m\|_{C^{2+\alpha}(S^{n-1})}\le K (\kappa K)^{m+1}\delta,
  $$
  $$
 \|P (\Phi_{m+1})-P(\Phi_m)\|_{C^{2+\alpha}(S^{n-1})}\le  (\kappa K)^{m+2}\delta.
 $$
 It follows that the sequence $\varphi_m$ converges in $C^{2+\alpha}(S^{n-1})$ to some function $\varphi\in C^{2+\alpha}(S^{n-1})$ 
 with
 $$
  \|\varphi\|_{C^{2+\alpha}(S^{n-1})}\le \frac{K\delta}{1-\kappa K}<\varepsilon
 $$
 if $\delta>0$ is small enough. This function  $\varphi$ will solve the equation $\widetilde{\Delta}\varphi=\gamma-P(\Phi)$, i.e., the function $f=1+\varphi$ will solve $Af=g$.
 
 We put $\varphi'=\varphi_0$ and $\varphi''=\varphi-\varphi_0$.
 It remains to estimate $\|\varphi''\|_{L^2(S^{n-1})}=$ $\|\varphi-\varphi_0\|_{L^2(S^{n-1})}$. To this end, we shall use \eqref{vc1} instead of \eqref{vc2} to obtain
 $$
 \|P(\Phi_0)\|_{L^2(S^{n-1})}=\|P(\Phi_0)-P(0)\|_{L^2(S^{n-1})}\le
 $$
 $$
 \le \kappa\max\limits_{i,j} \|(\Phi_0)_{x_ix_j}\|_{L^2(S^{n-1})}\le \kappa K\|\gamma\|_{L^2(S^{n-1})},
 $$
 so from the equation 
 $\widetilde{\Delta}(\varphi_1-\varphi_0)=-P(\Phi_0)$, we obtain
 $$
 \|\varphi_1-\varphi_0\|_{L^2(S^{n-1})}\le\|P(\Phi_0)\|_{L^2(S^{n-1})}\le \kappa K \|\gamma\|_{L^2(S^{n-1})}
 $$
 and
 $$
 \|(\Phi_1)_{x_ix_j}-(\Phi_0)_{x_ix_j}\|_{L^2(S^{n-1})}\le K\|P(\Phi_0)\|_{L^2(S^{n-1})}\le K(\kappa K)\|\gamma\|_{L^2(S^{n-1})}.
 $$
 Then
 $$
 \|P(\Phi_1)-P(\Phi_0)\|_{L^2(S^{n-1})}\le (\kappa K)^2 \|\gamma\|_{L^2(S^{n-1})},
 $$
 and we can continue as above to get inductively the inequalities
 $$
 \|\varphi_{m+1}-\varphi_m\|_{L^2(S^{n-1})}\le (\kappa K)^{m+1} \|\gamma\|_{L^2(S^{n-1})},
 $$
 $$
 \|(\Phi_{m+1})_{x_ix_j}-(\Phi_m)_{x_ix_j}\|_{L^2(S^{n-1})}\le K(\kappa K)^{m+1} \|\gamma\|_{L^2(S^{n-1})}
 $$
 (that requires the estimate $\max\limits_{i,j} \|(\Phi_m)_{x_ix_j}\|_{C(S^{n-1})}\le \frac{\kappa}{2C}$, but we have already obtained that bound even for the $C^{2+\alpha}$-norm of $\varphi_m$).
 
 Adding these estimates up, we get
 $$
  \|\varphi-\varphi_0\|_{L^2(S^{n-1})}\le\sum\limits_{m=0}^{\infty}\|\varphi_{m+1}-\varphi_m\|_{L^2(S^{n-1})}\le
 $$
 $$
 \sum\limits_{m=0}^{\infty} (\kappa K)^{m+1} \|\gamma\|_{L^2(S^{n-1})}=\frac{\kappa K}{1-\kappa K}\|\gamma\|_{L^2(S^{n-1})},
 $$
 and it remains to choose $\kappa>0$ so that 
 $\frac{\kappa K}{1-\kappa K}
 <\varepsilon$.
 \ep

\end{document}